\documentclass[12pt]{article}

\usepackage{amsmath,amssymb}

\newtheorem{theorem}{Theorem}

\newtheorem{fact}{Fact}
\newtheorem{proposition}{Proposition}
\newtheorem{lemma}{Lemma}

\newtheorem{remark}{Remark}
\newcommand{\proof}{{\bf Proof:~~~~}}
\newcommand{\qed}{ \hfill$\square$ }

\author{S. De Winter and M. Korb}

\title{Weak isometries of the Boolean cube}

\date{}

\begin{document}

\maketitle

\abstract{Consider the metric space $\mathcal{C}$ consisting of the $n$-dimensional Boolean cube equipped with the Hamming distance. A weak isometry of $\mathcal{C}$ is a permutation of $\mathcal{C}$ preserving a given subset of Hamming distances. In \cite{Krasin} Krasin showed that in most cases preserving a single Hamming distance forces a weak isometry to be an isometry.  In this article we  study those weak isometries that are not automatically an isometry, providing a complete classification of weak isometries of $\mathcal{C}$.}

\section{Introduction}

Distance preserving bijections of metric spaces, better known as {\it isometries}, play a prominent role throughout mathematics. Their weaker counterparts, not surprisingly called {\it weak isometries}, are bijections on a metric space that preserve only a prescribed subset of distances. This definition of weak isometry should not be confused with that of a weak isometry of a lattice ordered group as introduced by Jakub\'\i k \cite{Jakubik}.  In certain cases it is known that a weak isometry necessarily has to be an isometry. A famous example is the Beckman-Quarles Theorem that states that every mapping on the real Euclidean space $\mathbb{R}^n$, $n>1$, that preserves distance $1$ has to be an isometry \cite{BQ}. 

\medskip

In this article we will focus on weak isometries of the $n$-dimensional Boolean cube $\{0,1\}^n$ equipped with the Hamming distance $d$. We will denote this metric space by $\mathcal{C}_n$, or simply $\mathcal{C}$ if the dimension is understood.  Addition in this space is simply component wise addition modulo $2$. We recall that the Hamming distance between two elements (also called {\it words}) of $\mathcal{C}$ is simply the number of coordinates in which they differ.  We will use the notation $\overline{c}=(c_1,\hdots,c_n)$ for elements of $\mathcal{C}_n$,  $\overline{0}$ for the word $(0,\hdots,0)$, and $\overline{1}$ for the word $(1,\hdots,1)$. The {\it weight} of a word is simply its number of non-zero positions. Also, $E_n$ (or simply $E$ if the dimension $n$ is understood), respectively $O_n$ (or simply $O$ if the dimension $n$ is understood), will denote the set of all words of even, respectively odd, weight in $C_n$.

Let $\pi$ be any permutation of $\{1,2,\hdots,n\}$. The {\it permutation} of $\mathcal{C}_n$ induced by $\pi$ is the mapping on $\mathcal{C}_n$ defined by $ (x_1,\hdots,x_n)\mapsto (x_{\pi(1)},\hdots,x_ {\pi(n)})$. In order to simplify notation and terminology we will identify the permutation $\pi$ with the permutation of $\mathcal{C}_n$ it induces and write $\pi(x_1,\hdots,x_n)=(x_{\pi(1)},\hdots,x_ {\pi(n)})$. Clearly every permutation is an isometry.

Let $T_{\overline{a}}$ be the map $\mathcal{C}_n \rightarrow \mathcal{C}_n: \overline{c}\mapsto \overline{c}+\overline{a}$, where $\overline{a}$ is any fixed chosen word. Then $T_{\overline{a}}$ is obviously an isometry. We will call such maps {\it translations}.

We start by recalling the following well known fact (which also follows from our Lemma 1):

\begin{fact}
Let $\phi$ be an isometry of $\mathcal{C}_n$. Then there is a unique  permutation $\pi$  on $\{1,2,\hdots,n\}$, and a unique word $ \overline{a}$ of $\mathcal{C}_n$ such that 
\begin{equation}\label{isometry}
\phi=T_{\overline{a}}\circ \pi
\end{equation}
That is, every isometry is a composition of a translation and a permutation. Conversely, every mapping $\phi$ on $\mathcal{C}_n$ defined by (\ref{isometry}) above, where $\pi$ is any fixed permutation on $\{1,2,\hdots,n\}$, and where $T_{\overline{a}}$ is any translation, is an isometry of $\mathcal{C}_n$.
\end{fact}

A weak-isometry preserving all distances from a given set $P\subset\{1,2,\hdots,n\}$ will be denoted as a {\it $P$-isometry}. Whenever $P=\{p\}$ is a singleton we will simply write {\it $p$-isometry}. If $P$ consists of all even numbers from the set $\{1,2,\hdots,n\}$ we will call a $P$-isometry an {\it even isometry}.

In \cite{Krasin} Krasin proved the following theorem (an English translation of Krasin's original Russian paper appeared in \cite{KrasEng}):

\begin{theorem}[Krasin \cite{Krasin}]\label{Krasin}
If $\phi$ is a $p$-isometry of $\mathcal{C}_n$, with $p$ odd, and $p\notin\{\frac{n-1}{2},\frac{n}{2},\frac{n+1}{2},n\}$, then $\phi$ is a isometry.
\end{theorem}

The setup of his proof is very elegant, the main idea being that every such $p$-isometry must also be a $2$-isometry, which is proved using a nice counting argument. Two key steps to finish the prove are then that every $2$-isometry must be an even isometry and that every $1$-isometry must be an isometry. These two claims are rather obvious, but since no proof of them appears in Krasin's paper we will give proofs of them here that will be useful in the proof of our main result (see Lemma 1 and Lemma 2 of the next section). 

In Section \ref{Weak-isometries} we will look at $p$-isometries that do not satisfy the conditions of Theorem \ref{Krasin}. In \cite{Krasin} an example was given of such $p$-isometry that is not an isometry for each value of $p$. We will classify all weak-isometries. Our results can be summarized as follows: (here we refer to the theorems of section \ref{Weak-isometries} for the detailed classification of the various mappings)
\medskip

\noindent{\bf Main Theorem} {\it
Let $\phi$ be a weak-isometry of $\mathcal{C}_n$. Then $\phi$ is one of the following:
\begin{itemize}
\item an isometry;
\item an $\{n\}$-isometry (see Theorem \ref{niso});
\item an even isometry (see Theorem \ref{eveniso});
\item an $\{\frac{n}{2},n\}$-isometry (see Theorems \ref{half1} and \ref{half2});
\item an $\{\frac{n+1}{2}\}$-isometry (see Theorem \ref{nplusoneovertwo});
\item an $\{\frac{n-1}{2},\frac{n+1}{2},n\}$-isometry (see Theorem \ref{tripleiso}).
\end{itemize}}

\begin{remark} It is worthwhile to notice in the previous theorem that no $\frac{n-1}{2}$-isometries exist that are not $P$-isometries for some $P\supsetneq \{\frac{n-1}{2}\}$. This will be proved in subsection \ref{nminus1}
\end{remark}

\section{Two lemmas}

As mentioned in the introduction, the aim of this section is to give a proof of two simple lemmas on $p$-isometries with $p\in\{1,2\}$. The proof of the first lemma implicitly provides a proof of Theorem 1 (using graph theory Lemma 1 is immediate, but that approach does not provide the existence of the permutation $\pi$). The second lemma will also play an important role in the next section.

\begin{lemma}
Let $\phi$ be a $1$-isometry of $\mathcal{C}_n$, then $\phi$ is an isometry.
\end{lemma}
\proof We may assume without loss of generality that $\phi(\overline{0})=\overline{0}$. In case $\phi(\overline{0})=\overline{c}\neq\overline{0}$, just consider the $1$-isometry $T_{\overline{c}}\circ\phi$.

Hence, $\phi$ preserves all words of weight $1$, and as a consequence, the action of $\phi$ on the words of weight $1$ induces in a natural way a permutation $\pi$ on the coordinates. It is important to realize that at this point we only know that $\phi$ acts as $\pi$ on the words of weight $1$, and not necessarily on the words of higher weight. We will show that $\phi=\pi$ by proving that $\psi:=\pi^{-1}\circ\phi$ is the identity map.

Clearly $\psi$ fixes $\overline{0}$ and all word of weight $1$. We proceed by induction. Assume that $\psi$ fixes all words of weight less than or equal to $k-1$. Consider a word $\overline{c}$ of weight $k$.  Then there are exactly $k$ words of weight $k-1$ at distance $1$ from $\overline{c}$, and obviously $\psi(\overline{c})$ has to be a word of weight at least $k$ (as all words of weight less than $k$ are fixed) at distance $1$ from these $k$ words. The only possibility clearly is $\psi(\overline{c})=\overline{c}$. This proves the lemma. \qed

\begin{lemma}\label{even}
Let $\phi$ be a $2$-isometry of $\mathcal{C}_n$,  then $\phi$ is an even isometry.
\end{lemma}
\proof Without loss of generality we may assume that $\phi$ fixes $\overline{0}$, and hence permutes the words of weight $2$. If $n<4$ there is nothing to prove. If $n=4$ the word $\overline{1}$ is the unique word distinct from $\overline{0}$ that is at distance $2$ from all words of weight $2$. It hence has to be fixed and consequently $\phi$ is an even isometry. So assume that $n>4$.  We will start by constructing a permutation $\pi$ on the coordinates from this action on the words of weight $2$.

Consider all words of weight $2$ that have their first coordinate equal to $1$.  These words form a set of $n-1>3$ words of weight $2$ that are mutually at distance $2$. Now the maximal cardinality for a set of words of weight two that are mutually at distance $2$ but do not share a common coordinate equal to $1$ is easily seen to be $3$. Hence the images under $\phi$ of these $n-1$ words must have  a (unique) coordinate equal to $1$ in common, say coordinate $i$. Define $\pi(1)=i$. Repeating the argument for all coordinate positions produces a unique permutation $\pi$.  The action of $\phi$ and $\pi$ on the words of weight $2$ is clearly identical. We will now show that the action of $\phi$ and $\pi$ in fact coincide on all words of even weight. 

Let $E$ be the set of all words of even weight, and let $\phi_{E}$ and $\pi_E$ be the restriction of $\phi$ and $\pi$ to $E$. Note that this is well defined as $\phi$ and $\pi$  indeed induce permutations of the set $E$.
We will show that $\phi_E=\pi_E$ by proving that $\psi:=\pi_E^{-1}\circ\phi_E$ is the identity map on $E$.
Clearly $\psi$ fixes $\overline{0}$ and all word of weight $2$. We proceed by induction. Assume that $\psi$ fixes all words of even weight less than or equal to $2(k-1)$. Consider a word $\overline{c}$ of weight $2k$.  Then there are exactly $k(2k-1)$ words of weight $2(k-1)$ at distance $2$ from $\overline{c}$, and obviously $\psi(\overline{c})$ has to be a word of even weight at least $2k$ (as all words of even weight less than $2k$ are fixed) at distance $2$ from these $k(2k-1)$ words. The only possibility clearly is $\psi(\overline{c})=\overline{c}$.

Now let $\overline{c}$ be  any word of odd weight, and consider $T_{\phi(c)} \circ \phi \circ T_c$. This mapping is clearly a $2$-isometry fixing $\overline{0}$. Hence, by the above, we have that $T_{\phi(c)}\circ\phi\circ T_c=\sigma$ for some permutation $\sigma$ when restricted to $E$. It follows that $\phi T_c(\overline{w})=T_{\phi(c)}\sigma(\overline{w})$ for all words $\overline{w}\in E$. Hence, with $O$ being the set of all odd weight words, we see that $\phi T_c T_c(\overline{w})=T_{\phi(c)}\sigma T_c(\overline{w})$ for all words $\overline{w}\in O$. This implies that $\phi=T_{\phi(c)} \circ \sigma \circ T_c$ when restricted to $O$. Hence $\phi$ preserves all even distances, that is, $\phi$ is an even isometry. \qed

\section{Weak-isometries that are not isometries}\label{Weak-isometries}

We will now focus on permutations of $C$ that are weak-isometries but not isometries. Let $\epsilon$ be the set of non-zero even integers smaller than $n$. By Theorem \ref{Krasin} we know we can restrict ourselves to $P$-isometries where $P$ is a subset of $\{\frac{n-1}{2},\frac{n}{2},\frac{n+1}{2},n\}\cup\epsilon$.

We start with a simple yet useful lemma.

\begin{lemma}\label{useful}
\begin{itemize}
\item Every $\{p,n\}$-isometry is a $\{p,n-p,n\}$-isometry. 
\item Every $\{p,n-p\}$-isometry is a $\{p,n-p,n\}$-isometry.
\end{itemize}
\end{lemma}
\proof In order to prove the first statement, let $\overline{x}$ and $\overline{y}$ be two words at distance $n-p$, and let $\phi$ be a $\{p,n\}$-isometry. Then we see that $d(\overline{x},\overline{1}+\overline{y})=p$. Consequently $d(\phi(\overline{x}),\phi(\overline{1}+\overline{y}))=p$. Since $\phi$ is an $n$-isometry we also have that $\phi(\overline{1}+\overline{y})=\overline{1}+\phi(\overline{y})$, and hence $d(\phi(\overline{x}),\overline{1}+\phi(\overline{y}))=p$. It follows that $d(\phi(\overline{x}),\phi(\overline{y}))=n-p$, that is, $\phi$ preserves distance $n-p$.

Now let $\phi$ be a $\{p,n-p\}$-isometry, and let $\overline{x}$ and $\overline{y}$ be two words at distance $n$. Then one easily confirms that $\overline{y}$ is the unique word that is at distance $p$ from all words that are at distance $n-p$ from $\overline{x}$. It immediately follows that $\phi$ also has to preserve distance $n$. \qed

\subsection{The class of $n$-isometries}

\begin{theorem}\label{niso}
If $\phi$ is an $n$-isometry, then $\phi$ is a permutation of the pairs $\{\overline{c},\overline{1}+\overline{c}\}$, and conversely.
\end{theorem}
\proof This is obvious since for every word $\overline{c}$ there is a unique word at distance $n$, namely $\overline{1}+\overline{c}$, and the pairs $\{\overline{c},\overline{1}+\overline{c}\}$ partition $C$.\qed

\subsection{ The class of even isometries}

In \cite{Krasin} it was proven that every $p$-isometry, with $p$ even and $p\notin\{\frac{n-1}{2},\frac{n}{2},\frac{n+1}{2},n\}$ is a $2$-isometry. Invoking our Lemma \ref{even}  we obtain that every such $p$-isometry is in fact an even isometry. Furthermore the proof of Lemma \ref{even} provides the framework for the classification of all even isometries.

\begin{theorem}\label{eveniso}
If $\phi$ is an even isometry, then there exist two permutations, $\pi$ and $\sigma$, and two words $\overline{a}$ and $\overline{b}$, both of even weight, or both of odd weight, such that  
$$\phi=\left\{\begin{array}{ccc} T_{\overline{a}}\circ\pi & \mathrm{\ \ \ } & \mathrm{on\ } E\\   T_{\overline{b}}\circ  \sigma & \mathrm{\ \ \ } & \mathrm{on\ } O \end{array}\right., $$ and conversely.
\end{theorem}
\proof  Clearly there always is a unique word $\overline{c}$ (namely $\phi(\overline{0})$) such that $\phi=T_{\overline{c}}\circ \psi$, where $\psi$ is an even isometry fixing $\overline{0}$. By Lemma \ref{even} we know that 
$$\psi=\left\{\begin{array}{ccc} \pi & \mathrm{\ \ \ } & \mathrm{on\ } E\\   T_{\psi(\overline{d})} \circ \sigma \circ T_{\overline{d}} & \mathrm{\ \ \ } & \mathrm{on\ } O \end{array}\right.$$
where $\pi$ and $\sigma$ are two permutations, and $d$ is any word of odd weight. Hence
$$\phi=\left\{\begin{array}{ccc} T_{\overline{c}}\circ\pi & \mathrm{\ \ \ } & \mathrm{on\ } E\\   T_{\overline{c}}\circ T_{\psi(d)} \circ \sigma \circ T_d & \mathrm{\ \ \ } & \mathrm{on\ } O \end{array}\right.,$$ 
or, using $\sigma\circ T_{\overline{d}}=T_{\sigma(\overline{d})}\circ \sigma$,
$$\phi=\left\{\begin{array}{ccc} T_{\overline{c}}\circ\pi & \mathrm{\ \ \ } & \mathrm{on\ } E\\   T_{\overline{c}+ \psi(\overline{d}) + \sigma(\overline{d})}\circ \sigma & \mathrm{\ \ \ } & \mathrm{on\ } O \end{array}\right..$$
As there are no restrictions on which odd weight word $\psi(d)$ can be, as the sum of two odd weight words is always even, and as there are no restrictions on which word $\phi({\overline{0}})$ can be, we obtain that an even isometry $\phi$ is always of the form 
$$\phi=\left\{\begin{array}{ccc} T_{\overline{a}}\circ\pi & \mathrm{\ \ \ } & \mathrm{on\ } E\\   T_{\overline{b}}\circ  \sigma & \mathrm{\ \ \ } & \mathrm{on\ } O \end{array}\right.$$
where $a$ and $b$ are  any two words both of even weight or both of odd weight, and $\pi$ and $\sigma$ are any two permutations. Conversely it is obvious that any such mapping is an even isometry.\qed

\subsection{The class of $\{\frac{n}{2},n\}$-isometries}

By Lemma \ref{useful} we know that every $\frac{n}{2}$-isometry automatically is a $\{\frac{n}{2},n\}$-isometry.
Now let $\phi$ be a $\{\frac{n}{2},n\}$-isometry fixing $\overline{0}$ (this is not a real restriction since every $\{\frac{n}{2},n\}$-isometry is like this up to translation). Then, as $\phi$ is an $n$-isometry, $\phi$ has to permute the pairs $\{\overline{c},\overline{1}+\overline{c}\}$. Now each such pair contains a word of weight less than or equal to $n/2$. Hence, if we identify each of the pairs that contain a word of weight less than $n/2$ with that word, then $\phi$ naturally induces a mapping on the words of weight less than $n/2$. Also, the pairs $\{\overline{c},\overline{1}+\overline{c}\}$ consisting of two words of weight $n/2$ contain a unique word that has a zero entry in its first coordinate. Hence, if we identify each of the pairs that contain a word of weight  $n/2$ with that word, then $\phi$ naturally induces a mapping on the words of weight $n/2$ with first coordinate equal to zero. Define $X:=X_1\cup X_2$ where $X_1$ is the set of all words of weight less than $n/2$ and $X_2$ is the set of all words of weight $n/2$ that have their first coordinate equal to zero. We will now concentrate on the map $\psi$ which acts on $X$ through the induced action of $\phi$ as described above.  Clearly $\psi$ is a permutation of $X$. Furthermore $\psi$ fixes $\overline{0}$, and preserves distance $n/2$, since if $d(\overline{x},\overline{y})=n/2$ then also $d(\overline{1}+\overline{x},\overline{y})=d(\overline{x},\overline{1}+\overline{y})=d(\overline{1}+\overline{x},\overline{1}+\overline{y})=n/2$.   
\medskip

\noindent{\bf CASE I} We first deal with the case where $n/2$ is odd. Set $A_{2k}$ equal to the number of words of weight $n/2$   at distance $n/2$  from a fixed word of weight $0<2k<n/2$. We deduce (see also \cite{Krasin}) that
$$A_{2k}= {2k\choose k}{n - 2k \choose \frac{n}{2} -k}. $$ 
The proof of Proposition 3 of \cite{Krasin} now yields $A_2>A_4>\hdots>A_{\frac{n}{2}-1}$. As a consequence $\psi$ has to preserve the weight of all even weight words in $X$. Now define $B_{2k+1}$, $0<2k+1<n/2$, to be the number of possible weights for a word in $X$ at distance $n/2$ from a word of weight $2k+1$. One easily sees that $B_{2k+1}=k+1$. Furthermore, since a word at distance $n/2$ from a word of weight $2k+1$ must have even weight, and $\psi$ does not change the weight of even weight words, we obtain, using $B_{2k+1}\neq B_{2l+1}$ if $k\neq l$, that $\psi$ is weight preserving. Hence from the action of $\psi$ on the words of weight $1$ we obtain in a natural way a permutation $\pi$.  We now proceed by proving that $\psi$ acts as $\pi$ on $X$. To this end it is sufficient to proof that if $\pi$ is the identity, then so is $\psi$ (where some care is required for the action on the words of weight $n/2$, see below). So from here we suppose that $\psi$ fixes all words of weight $1$.
Consider any word $\overline{c} $ of weight $n/2-1$, and the $n/2+1$ words of weight $1$ whose support is not contained in the support of $\overline{c}$. Then $\overline{c}$ is the unique word of weight $n/2-1$ at distance $n/2$ from all of  the $n/2+1$ considered words of weight $1$. Hence $\overline{c}$ has to be fixed, and consequently all words of weight $n/2-1$ have to be fixed. Now let $\overline{c}$ be any word of weight $1<2k+1< n/2$, and let $D(\overline{c})$ be the set of words of weight $n/2-1$ at distance $n/2$ from $\overline{c}$. Then $D(\overline{c})$ is the set of all words of weight $n/2-1$ that have exactly $k$ non-zero positions in the support of $\overline{c}$. One easily deduces that  $D(\overline{x})\neq D(\overline{y})$ if $\overline{x}\neq\overline{y}$ (note that $2k+1<n/2$ is necessary here). Hence $\overline{c}$ has to be fixed, and consequently all words of odd weight less than $n/2-1$ in $X$ are fixed. Specifically all words of weight $n/2-2$ are fixed. Now let $\overline{c}$ be any word of weight $1<2k< n/2-1$, and define $F(\overline{c})$ to be the set of words of weight $n/2-2$ at distance $n/2$ from $\overline{c}$. Then $F(\overline{c})$ is the set of all words of weight $n/2-2$ that have exactly $k-1$ non-zero positions in the support of $\overline{c}$. One easily deduces that  $F(\overline{x})\neq F(\overline{y})$ if $\overline{x}\neq\overline{y}$. Hence $\overline{c}$ has to be fixed, and consequently all words of even weight in $X$ are fixed. Finally let $\overline{c}$ be any word of weight $n/2$. Then $D(\overline{c})= D(\overline{y})$ if and only if $\overline{y}\in\{\overline{c},\overline{1}+\overline{c}\}$. Hence $\overline{c}$ is fixed. (Note that if $\pi$ is not the identity, then one obtains the conclusion $\psi(\overline{c})=\pi(\overline{c})$ if $(\pi(\overline{c}))_1=0$, and $\psi(\overline{c})=\pi(\overline{c})+\overline{1}$ if $(\pi(\overline{c}))_1=1$.) Hence we completely understand the action of $\psi$ on $X$. From this we now easily deduce the following

\begin{theorem}\label{half1}
Suppose $2\mid n$ but $4\not\:\mid n$. If $\phi$ is a $\{n/2,n\}$-isometry fixing $\{\overline{0},\overline{1}\}$, then there is a permutation $\pi$, and a subset $S$ of $X$ such that 
\begin{eqnarray}\label{type1}\phi(\overline{c})=\left\{\begin{array}{ccc} \pi(\overline{c}) & \mathrm{\ \ \ } & \mathrm{if \  } \{\overline{c},\overline{1}+\overline{c}\}\cap S=\emptyset\\ 
\pi(\overline{c})+\overline{1} & \mathrm{\ \ \ } & \mathrm{ if \ } \{\overline{c},\overline{1}+\overline{c}\}\cap S\neq\emptyset    \end{array}\right..\end{eqnarray}
Conversely, for any permutation $\pi$ and any $S\subset X$, Equation \ref{type1} defines an $\{n/2,n\}$-isometry fixing $\{\overline{0},\overline{1}\}$.
\end{theorem}

\noindent{\bf CASE II} Next we consider the case where $n/2$ is even. The major issue here is that a word of odd weight can never be at distance $n/2$ from a word of even weight. In exactly the same way as in CASE I we obtain that $\psi$ has to preserve the weight of all even weight words in $X$. We will now concentrate on the action of $\psi$ on $E_X$, the subset of $X$ consisting of even weight words. Consider two words of weight two that share one position in their support. Then there are exactly $${{n-3}\choose{\frac{n}{2}-1}} + {{n-3}\choose{\frac{n}{2}-2}}$$ words of weight $n/2$ at distance $n/2$ from both of these words. Now consider two words of weight two with disjoint supports. Then there are exactly $$4{{n-4}\choose{\frac{n}{2}-2}}$$ words of weight $n/2$ at distance $n/2$ from both of these words. One easily checks that 
$${{n-3}\choose{\frac{n}{2}-1}} + {{n-3}\choose{\frac{n}{2}-2}}=4{{n-4}\choose{\frac{n}{2}-2}}\frac{n-3}{n-2}$$ from which follows that $\psi$ has to map words of weight two with intersecting support to words of weight two with intersecting support. We deduce, in the same manner as in Lemma \ref{even}, that the action of $\psi$ on the words of weight two coincides with the action of a permutation, say $\pi$. We will now show that the action of $\psi$ on $E_X$ coincides with $\pi$ (where again some care is required for the action on the words of weight $n/2$). To this end it is sufficient to proof that if $\pi$ is the identity, then so is $\psi$. So from here we suppose that $\psi$ fixes all words of weight $2$. Consider all words of weight two that have their support in a given set of $n/2+2$ positions. Then there is a unique word of weight $n/2-2$ at distance $n/2$ from all these words. It follows that all words of weight $n/2-2$ are fixed, and hence, using the same type of argument as under CASE I, that all words of even weight less than $n/2$ in $X$ are fixed. Finally consider a word $\overline{c}$ of weight $n/2$ in $X$. Define $F(\overline{c})$ to be the set of words of weight $n/2-2$ at distance $n/2$ from $\overline{c}$. Then$F(\overline{c})$ is the set of words of weight $n/2-2$ whose support shares exactly $n/4-1$ positions with the support of $\overline{c}$. Then $F(\overline{c})= F(\overline{y})$ if and only if $\overline{y}\in\{\overline{c},\overline{1}+\overline{c}\}$. Hence $\overline{c}$ is fixed, and we completely understand the action of $\psi$ on $E_X$. 
Since $n/2$ and $n$ are even, the action of $\phi$ on $E$ and on $O$ are completely independent. Furthermore, $E=O+\overline{c}$ for every odd weight word $\overline{c}$. We obtain

\begin{theorem}\label{half2}
Suppose $4\mid n$. If $\phi$ is a $\{n/2,n\}$-isometry fixing $\{\overline{0},\overline{1}\}$, then there are two permutation $\pi_1$ and $\pi_2$, two subsets $S_1$ and $S_2$ of $E_X$, and two odd weight words $\overline{a}$ and $\overline{b}$ such that 
\begin{eqnarray}\label{type2}\phi(\overline{c})=\left\{\begin{array}{ccc} \pi_1(\overline{c}) & \mathrm{\ \ \ } & \mathrm{if \  } \{\overline{c},\overline{1}+\overline{c}\}\cap S_1=\emptyset\\ 
\pi_1(\overline{c})+\overline{1} & \mathrm{\ \ \ } & \mathrm{ if \ } \{\overline{c},\overline{1}+\overline{c}\}\cap S_1\neq\emptyset \\ \overline{b}+\pi_2(\overline{c}+\overline{a}) & \mathrm{\ \ \ } & \mathrm{if \  } \{\overline{c}+\overline{a},\overline{1}+\overline{c}+\overline{a}\}\cap S_2=\emptyset\\
  \overline{b}+\pi_2(\overline{c}+\overline{a})+\overline{1} & \mathrm{\ \ \ } & \mathrm{ if \ } \{\overline{c}+\overline{a},\overline{1}+\overline{c}+\overline{a}\}\cap S_2\neq\emptyset  \end{array}\right..\end{eqnarray}
Conversely, for any two permutation $\pi_1$ and $\pi_2$, any two $S_1,S_2\subset E_X$, and any two odd weight words $\overline{a}$ and $\overline{b}$, Equation \ref{type2} defines an $\{n/2,n\}$-isometry fixing $\{\overline{0},\overline{1}\}$.
\end{theorem}

This classifies all $\{n/2,n\}$-isometries.

\subsection{The class of $\{\frac{n+1}{2}\}$-isometries}

In this section we will show that all $\{\frac{n+1}{2}\}$-isometries were essentially determined in the previous section.
Let $\phi$ be an $\{\frac{n+1}{2}\}$-isometry on $\mathcal{C}_n$. Embed $\mathcal{C}_n$  in $\mathcal{C}_{n+1}$ as the subspace that has its first coordinate equal to $0$, that is, a word $\overline{c}=(c_1,c_2,\hdots,c_n)$ of $\mathcal{C}_n$ is embedded in $\mathcal{C}_{n+1}$ as the word $(0,\overline{c})=(0,c_1,c_2,\hdots,c_n)$. Now define $\psi$ to be the following map on $\mathcal{C}_{n+1}$:

$$\left\{\begin{array}{ccc} \psi((0,\overline{c})) & = & (0,\phi(\overline{c})) \\ \psi(\overline{1}+(0,\overline{c})) & = & \overline{1}+(0,\phi(\overline{c}))\end{array}\right.$$

Then $\psi$ acts as $\phi$ on the embedded sunspace $\mathcal{C}_n$, and furthermore is an $(n+1)$-isometry.

\begin{lemma}
The map $\psi$ is an $\{\frac{n+1}{2},n+1\}$-isometry on $\mathcal{C}_{n+1}$.
\end{lemma}
\proof In view of  the remark preceeding the lemma we only need to prove that $\psi$ is an $\frac{n+1}{2}$-isometry. To this end let $\overline{x}$ and $\overline{y}$ be two words of $\mathcal{C}_{n+1}$ at distance $\frac{n+1}{2}$. If both words belong to the embedded subspace $\mathcal{C}_n$ then $d(\psi(\overline{x}),\psi(\overline{y}))=\frac{n+1}{2}$ follows from the fact that $\phi$ is an $\frac{n+1}{2}$-isometry on $\mathcal{C}_n$. Next assume that both words belong to the complement of $\mathcal{C}_n$ in $\mathcal{C}_{n+1}$. Using $\overline{1}+\psi(\overline{x})=\psi(\overline{1}+\overline{x})$, we obtain 
$$\begin{array}{ccl} 

\frac{n+1}{2} & = & d(\overline{x},\overline{y})\\  & = & d(\overline{1}+\overline{x},\overline{1}+\overline{y})\\  & = & d(\psi(\overline{1}+\overline{x}),\psi(\overline{1}+\overline{y})) \\  & = & d(\overline{1}+\psi(\overline{x}),\overline{1}+\psi(\overline{y}))   \\  & = &  d(\psi(\overline{x}),\psi(\overline{y})).\end{array}$$

Finally assume that $\overline{x}$ belongs to the embedded subspace $\mathcal{C}_n$, and $\overline{y}$ belongs to its complement. Using the above we obtain
$$\begin{array}{ccl} d(\psi(\overline{x}),\psi(\overline{y})) & = & d(\psi(\overline{x}),\overline{1}+\psi(\overline{1}+\overline{y})) \\  & = & n+1- d(\psi(\overline{x}),\psi(\overline{1}+\overline{y})) \\  & = & n+1- d(\overline{x},\overline{1}+\overline{y}) \\  & = &n+1- \frac{n+1}{2}  \\& = & \frac{n+1}{2}.\end{array}.$$
This proves the lemma.
\qed
\medskip

Conversely it is obvious that any $\{\frac{n+1}{2},n+1\}$-isometry on $\mathcal{C}_{n+1}$ that stabilizes the embedded $\mathcal{C}_n$ induces an $\frac{n+1}{2}$-isometry on $\mathcal{C}_n$. So in order to determine all $\frac{n+1}{2}$-isometries of $\mathcal{C}_n$ we only need to determine which $\{\frac{n+1}{2},n+1\}$-isometries on $\mathcal{C}_{n+1}$ stabilize $\mathcal{C}_n$. 
\medskip

In order to present our theorem as concise as possible we first introduce the concept of a {\it $\sigma_{i,j}$-mapping}.  The proof of the next theorem will clarify why these mappings arise naturally. They will also play an important role in the classification of $\{\frac{n-1}{2},\frac{n+1}{2},n\}$-isometries.
\medskip

\noindent {\bf Definition.}  Let $n$ be fixed, and let $i,j\in\{1,2,\hdots,n\}$. Let $\sigma$ be a bijection from $ \{1,2,\hdots,n\}\setminus\{i\}$ to  $\{1,2,\hdots,n\}\setminus\{j\}$. Then $\phi$ is a {\it $\sigma_{i,j}$ mapping on $\mathcal{C}_n$} if and only if 
$$\phi: \mathcal{C}_n\rightarrow\mathcal{C}_n, \  \overline{c}=(c_1,c_2,\hdots,c_n)\mapsto \phi(\overline{c})=(d_1,d_2,\hdots,d_n),$$
where
\begin{itemize}
\item if $c_i=0$, then $d_k=c_{\sigma^{-1}(k)}$ for $k\neq j$, and $d_j=0$;
\item if $c_i=1$, then $d_k=c_{\sigma^{-1}(k)}+1$ for $k\neq j$, and $d_j=1$.
\end{itemize}

\medskip

Note that we do not require that $i\neq j$.
The following proposition is a direct consequence of the above definition:

\begin{proposition}
Let $\phi$ be a $\sigma_{ij}$-mapping, then
\begin{itemize}
\item $\phi(E)=E$ and $\phi(O)=O$;
\item $\phi\circ T_{\overline{a}}=T_{\phi(\overline{a})}\circ\phi$.
\end{itemize}
\end{proposition}

In view of the previous section we obtain the following.

\begin{theorem}\label{nplusoneovertwo}
If $n\equiv 1 \pmod{4}$ and $\phi$ is an $\frac{n+1}{2}$-isometry of $\mathcal{C}_n$, fixing $\overline{0}$, but such that $\phi$ is not an isometry,  then $\phi$ is a $\sigma_{i,j}$-mapping on $\mathcal{C}_n$. Conversely every   $\sigma_{i,j}$-mapping on $\mathcal{C}_n$ is an $\frac{n+1}{2}$-isometry.
\medskip

\noindent If $n\equiv 3 \pmod{4}$ and $\phi$ is an $\frac{n+1}{2}$-isometry of $\mathcal{C}_n$, fixing $\overline{0}$, but such that $\phi$ is not an even isometry, then there are three possibilities:
\begin{itemize}

\item[I)] there is a $\sigma_{i,j}$-mapping $\tau$ on $\mathcal{C}_n$ and an even weight word $\overline{a}$, such that $\phi$ acts as an even isometry  on the words of even weight and as $T_{\overline{a}}\circ\tau$ on the words of odd weight;

\item[II)]  there is a $\sigma_{i,j}$-mapping $\tau$ on $\mathcal{C}_n$, such that $\phi$ acts as $\tau$ on the words of even weight and as an even isometry  on the words of odd weight;

\item[III)] there is a $\sigma_{i,j}$-mapping $\tau_1$ on $\mathcal{C}_n$, a   $\sigma_{k,l}$-mapping $\tau_2$ on $\mathcal{C}_n$ and an even weight word $\overline{a}$, such that $\phi$ acts as $\tau_1$ on the words of even weight, and as $T_{\overline{a}}\circ\tau_2$ on the words of odd weight.

\end{itemize}

\end{theorem}
\proof Let  $n\equiv 1 \pmod{4}$ and let  $\phi$ be an $\frac{n+1}{2}$-isometry of $\mathcal{C}_n$ fixing $\overline{0}$. Embed $\mathcal{C}_n$ in $\mathcal{C}_{n+1}$ and define $\psi$ as earlier in this section. Then, invoking Theorem \ref{half1}, we need to determine the restriction on $\pi$ and $S$ in order for $\psi$ to stabilize $\mathcal{C}_n$ in $\mathcal{C}_{n+1}$. We easily see that there are no restriction on $\pi$, and that, if $\pi(0)=0$ then $S$ must be the empty set and $\phi$ is a usual isometry. If $\pi^{-1}(0)=i\neq0$, then $c_i=1$ must imply that $\{\overline{c},\overline{c}+\overline{1}\}\cap S\neq\emptyset$, whereas $c_i=0$ must imply that $\{\overline{c},\overline{c}+\overline{1}\}\cap S=\emptyset$ (here the word $\overline{c}$ is any word of $\mathcal{C}_{n+1}$ with first coordinate zero). Given $\pi$ this defines $S$ uniquely. Formulating this directly in $\mathcal{C}_n$ yields a $\sigma_{i,j}$ map (with $j=\pi(0)$), and hence the theorem.

The case $n\equiv 3 \pmod{4}$ is analogues. Here one invokes Theorem \ref{half2} and recalls that the action on the even and odd weight words is totally independend. When the action on the odd weight words is not an even isometry, we obtain that there must be two odd weight words $\overline{a}$ and $\overline{b}$, and a $\sigma_{ij}$-mapping $\tau$ such that $\phi$ acts on $O$ as $T_{\overline{b}}\circ \tau \circ T_{\overline{a}}$. In view of the proposition preceding the theorem this is equivalent with the existence of an even weight word $\overline{c}$ ($=\overline{b}+\tau(\overline{a})$) such that $\phi$ acts like $T_{\overline{c}}\circ\tau$ on $O$. This proves the theorem.\qed

\subsection{The class of $\{\frac{n-1}{2}\}$-isometries}\label{nminus1}

It is the aim of this section to show that every $\{\frac{n-1}{2}\}$-isometry is in fact a $P$-isometry for some $P\supsetneq \{\frac{n-1}{2}\}$. This will be obtained by extending an argument that was made in \cite{Krasin}.
Following Krasin \cite{Krasin}, let $A_{2k}$ equal the number of words of weight $\frac{n-1}{2}$ at distance $\frac{n-1}{2}$ from a given word of weight $2k$. We easily compute  that $$A_{2k}={2k \choose k} {n-2k\choose \frac{n-1}{2} -k}.$$
Then, still using the notation of \cite{Krasin}, $$H_k=\frac{A_{2(k+1)}}{A_{2k}}=\frac{(2k+1)(\frac{n-1}{2}-k+1)}{(k+1)(n-2k)}.$$
A simple computation now shows that $H_k<1$ if and only if $2k<\frac{n-1}{2}$ and $H_k>1$ if and only if $2k>\frac{n-1}{2}$, and hence $$A_2>A_4>\cdots>A_{2t} \hspace*{1cm}\mathrm{and}\hspace*{1cm} A_{2t+2}<A_{2t+4}<\cdots<A_{n-1},$$ where $t=\frac{n-1}{4}$ if $n\equiv 1 \pmod{4}$, and $t=\frac{n+1}{4}$ if $n\equiv 3 \pmod{4}$.
\medskip

Now let $\phi$ be an $\{\frac{n-1}{2}\}$-isometry fixing $\overline{0}$. Then clearly $\phi$ can map a word of weight $2k$ to a word of weight $2l$ if and only if $A_{2k}=A_{2l}$.

\begin{theorem}
If $n\equiv 3 \pmod{4}$ and $\phi$ is an $\frac{n-1}{2}$-isometry then $\phi$ is an $\{\frac{n-1}{2},\frac{n+1}{2},n\}$-isometry.

If $n\equiv 1 \pmod{4}$ and $\phi$ is an $\frac{n-1}{2}$-isometry then $\phi$ is an even isometry.
\end{theorem}
\proof First consider the case $n\equiv 3 \pmod{4}$, and  let $\phi$ be an $\{\frac{n-1}{2}\}$-isometry fixing $\overline{0}$ (we can assume the latter without loss of generality).  Now a simple direct computation shows that $A_{2k}=A_{n-2k+1}$. It follows that $A_{2l}=A_{\frac{n+1}{2}}$ implies that $l=\frac{n+1}{4}$, and hence that $\phi$ maps words of weight $\frac{n+1}{2}$ to words of weight $\frac{n+1}{2}$. Now let $\overline{x}$ and $\overline{y}$ be any two words at distance $\frac{n+1}{2}$. Define $$\psi:=T_{\phi(\overline{x})}\circ\phi \circ T_{\overline{x}}.$$ Clearly $$\phi= T_{\phi(\overline{x})}\circ\psi \circ T_{\overline{x}},$$ and consequently $$d(\phi(\overline{x}),\phi(\overline{y}))=d(\psi(\overline{x}),\psi(\overline{y})).$$  Now $\phi \circ T_{\overline{x}}$ maps $\overline{x}$ to $\overline{0}$, and $\overline{y}$ to a word of weight $\frac{n+1}{2}$ by the above. Hence $$d(\psi(\overline{x}),\psi(\overline{y}))=\frac{n+1}{2}$$ and we conclude that $\phi$ is an $\frac{n+1}{2}$-isometry. Using Lemma \ref{useful} we obtain that $\phi$ is an  $\{\frac{n-1}{2},\frac{n+1}{2},n\}$-isometry.
\medskip

Next consider the case $n\equiv 1 \pmod{4}$. Let $\phi$ be an $\{\frac{n-1}{2}\}$-isometry fixing $\overline{0}$ (we can assume the latter without loss of generality). As before $A_{2k}=A_{n-2k+1}$. We obtain that $A_{2l}=A_{\frac{n+3}{2}}$ implies that $l\in\{\frac{n-1}{2},\frac{n+3}{2}\}$. However, since $\phi$ permutes the words of weight $\frac{n-1}{2}$ it follows that $\phi$ also must permute the words of weight $\frac{n+3}{2}$. As in the first part of the proof we can deduce from this that $\phi$ is an $\frac{n+3}{2}$-isometry. It follows that $\phi$ is an even isometry. This proofs the theorem.\qed

\subsection{The class of $\{\frac{n-1}{2},\frac{n+1}{2},n\}$-isometries}

Clearly every $\{\frac{n-1}{2},\frac{n+1}{2},n\}$-isometry is also a $\frac{n+1}{2}$-isometry. Hence we will obtain a classification of $\{\frac{n-1}{2},\frac{n+1}{2},n\}$-isometries by finding which of the mappings described in Theorem \ref{nplusoneovertwo} are also $n$-isometries.

\begin{theorem}\label{tripleiso}
Let $\phi$ be an $\{\frac{n-1}{2},\frac{n+1}{2},n\}$-isometry fixing $\overline{0}$.

If $n\equiv 1 \pmod{4}$ then $\phi$ is an isometry.

If $n\equiv 3 \pmod{4}$ then $$\phi=\left\{\begin{array}{ccc}\tau & & \mathrm{on\  }E\\ T_{\overline{1}^j}\circ\tau & & \mathrm{on\ }O \end{array}\right.,$$
where $\overline{1}^j$ is the word with $1$ in all positions except for a unique $0$ in position $j$, and where $\tau$ is a $\sigma_{ij}$-mapping.
\end{theorem}
\proof Let $n\equiv 1 \pmod{4}$, and assume that $\phi$ is not an isometry. Let $i,j$ and $\sigma$ be as in Theorem \ref{nplusoneovertwo}. Consider a word $\overline{c}=(c_1,c_2,\hdots,c_n)$ with $c_i=0$. Set $\overline{d}=\phi(\overline{c})$. Then $d_k=c_{\sigma^{-1}(k)}$ for $k\neq j$. Now the unique word at distance $n$ from $\overline{c}$ clearly has its $i$th coordinate equal to $1$.  Set $\overline{e}=\phi(\overline{c}+\overline{1})$. Then  $e_k=1+(c_{\sigma^{-1}(k)}+1)=c_{\sigma^{-1}(k)}$ for $k\neq j$, implying that $\overline{e}\neq\overline{d}+\overline{1}$, and hence that $\phi$ is not an $n$-isometry. This proofs the first part of the theorem.

Now let $n\equiv 3 \pmod{4}$,  and assume that $\phi$ is not an even isometry. We need to check which of the three types of mappings  described in the second part of Theorem \ref{nplusoneovertwo} might give rise to an $\{\frac{n-1}{2},\frac{n+1}{2},n\}$-isometry fixing $\overline{0}$. However, it is easily seen that any map that acts as an even isometry on either the even or odd weight words and also has to be an $n$-isometry has to be an even isometry on the whole space ($n$ is odd!). Hence, only the mappings described in case III) of the second part of Theorem \ref{nplusoneovertwo} might give rise to an $\{\frac{n-1}{2},\frac{n+1}{2},n\}$-isometry fixing $\overline{0}$.  Let $\phi$ be such mapping, and let $\overline{c}$ and $\overline{c}+\overline{1}$ be two words at distance $n$. Without loss of generality assume that $\overline{c}$ has even weight. Then, by our assumptions, $\phi(\overline{c})=\tau_1(\overline{c})$ and $\phi(\overline{c}+\overline{1})=T_{\overline{a}}\tau_2 (\overline{c}+\overline{1})$,  for some $\sigma_{i,j}$-mapping $\tau_1$,  some $\sigma_{k,l}$-mapping $\tau_2$, and fixed even weight word $\overline{a}$. For simplicity of notation we will also use $\tau_1$ and $\tau_2$ for the bijections on the coordinate positions corresponding with these mappings. We need to check under which conditions $\phi(\overline{c})+\overline{1}=\phi(\overline{c}+\overline{1})$, that is, under which conditions $\tau_1(\overline{c})+\overline{1}=T_{\overline{a}}\tau_2(\overline{c}+\overline{1})$.  Suppose that $c_i=0$. Then $(\tau_1(\overline{c}))_t=c_{\tau_1^{-1}(t)}$ for $t\neq j$ and $(\tau_1(\overline{c}))_j=0$.
On the other hand we have
\begin{eqnarray}\label{coef1}(T_{\overline{a}}\tau_2(\overline{c}+\overline{1}))_t= a_t+(\overline{c}+\overline{1})_{\tau_2^{-1}(t)}=a_t+c_{\tau_2^{-1}(t)}+1,\end{eqnarray}
 or 
\begin{eqnarray}\label{coef2}(T_{\overline{a}}\tau_2(\overline{c}+\overline{1}))_t= a_t+(\overline{c}+\overline{1})_{\tau_2^{-1}(t)}+1=a_t+c_{\tau_2^{-1}(t)}\end{eqnarray}
 for $t\neq l$, depending on whether $(\overline{c}+\overline{1})_k=c_k+1$ equals $0$ or $1$. We will show that $i=k$. Assume that $i\neq k$. Then we can always find a position $t\neq l$, and two even weight words $\overline{c}$ and $\overline{d}$ such that $c_i=d_i=0$, $c_k\neq d_k$ and $c_{\tau_2^{-1}(t)}=d_{\tau_2^{-1}(t)}$. Using Equations (\ref{coef1}) and (\ref{coef2}) we then see that $\phi(\overline{c}+\overline{1})_t\neq\phi(\overline{d}+\overline{1})_t$, and hence that $\phi$ cannot be an $n$-isometry. If follows that $i=k$.
\medskip

Next we need to assure that $(T_{\overline{a}}\tau_2(\overline{c}+\overline{1}))_j=1$ independent of our choice of $\overline{c}$. Assume that $j\neq l$. Then, using Equation (\ref{coef1}) and (\ref{coef2}), we see that $(T_{\overline{a}}\tau_2(\overline{c}+\overline{1}))_j$ depends on $c_{\tau_2^{-1}(j)}$, which depends on our choice of $\overline{c}$ since $\tau_2^{-1}(j)\neq i$ (recall that because $i=k$, $i$ is not a preimage for $\tau_2$). We conclude that $j=l$, and hence both $\tau_1$ and $\tau_2$ are $\sigma_{ij}$-mappings.
\medskip

We will show that in fact $\tau_1=\tau_2$. Since $i=k$ we have, for $t\neq j$, $$(T_{\overline{a}}\tau_2(\overline{c}+\overline{1}))_t= a_t+(\overline{c}+\overline{1})_{\tau_2^{-1}(t)}+1=a_t+c_{\tau_2^{-1}(t)}$$ which should equal $(\tau_1(\overline{c}))_t+1=c_{\tau_1^{-1}(t)}+1$, or equivalently,
\begin{eqnarray}\label{coef3} a_t+c_{\tau_2^{-1}(t)}=c_{\tau_1^{-1}(t)}+1,\end{eqnarray}
for $t\neq j$ and any $\overline{c}$ with $c_i=0$. If $\tau_1$ were not equal to $\tau_2$ then there is a $t$ such that $\tau_1^{-1}(t)\neq\tau_2^{-1}(t)$, say $\tau_1(u)=\tau_2(v)=t$, with $u\neq v$. Since $n\geq7$ (note that the case $n=3$ is trivial), we can find two words of even weight $\overline{c}$ and $\overline{d}$ with $c_i=d_i=0$, $c_u=c_v=d_u\neq d_v$. But then Equation (\ref{coef3}) gives 
$$a_t+c_v=c_u+1=d_u+1=a_t+d_v,$$ a contradiction since $c_v\neq d_v$. We conclude that $\tau_1=\tau_2$. For simplicity of notation set $\tau_1=\tau$. Now Equation(\ref{coef3}) implies that $a_t=1$ for all $t\neq j$. Finally, $\tau(\overline{c})_j+1$ equals $1$, and hence so should $(T_{\overline{a}}\tau (\overline{c}+\overline{1}))_j$. Using $(\overline{c}+\overline{1})_i=1$ this easily yields $a_j=0$. Consequently $\overline{a}=\overline{1}^j$, that is, the word with $1$ in all positions except for a unique $0$ in position $j$. 

\medskip

Finally we need to show that $$\phi=\left\{\begin{array}{ccc}\tau & & \mathrm{on\  }E\\ T_{\overline{1}^j}\circ\tau & & \mathrm{on\ }O \end{array}\right.$$ is indeed an 
$\{\frac{n-1}{2},\frac{n+1}{2},n\}$-isometry for every $\sigma_{ij}$-mapping $\tau$. This map is clearly an $\frac{n+1}{2}$-isometry, so it is sufficient to show that it is also an $n$-isometry, that is, we need to show that $\phi(\overline{c})+\overline{1}=\phi(\overline{c}+\overline{1})$. Without loss of generality we can assume that $\overline{c}\in E$. There are two cases to consider, namely $c_i=0$ and $c_i=1$. However in both cases a straightforward computation shows that indeed $\phi(\overline{c})+\overline{1}=\phi(\overline{c}+\overline{1})$. This proves the theorem. \qed
\bigskip

\begin{remark}
It is worthwhile to note that the map $I_i^-$ given by Krasin as an example of an $\frac{n+1}{2}$-isometry is the above map with $i=j$ and with as bijection underlying the map $\tau$ the identity.
\end{remark}

Putting all results of this section together now gives the Main Theorem.

\end{document}